\documentclass[a4paper,leqno,12pt]{article}

\usepackage[pdftex]{graphics,color}
\usepackage{graphics}
\usepackage{amsmath,bm,stmaryrd,mathtools, xy, subfigure}
\usepackage{array}
\usepackage{amscd}
\usepackage{amssymb}
\usepackage{enumerate}
\usepackage{indentfirst}
\usepackage{latexsym}
\usepackage{multicol}
\usepackage{color}
\usepackage{theorem}

\makeatletter\@addtoreset{equation}{section}\makeatother

\newtheorem{theorem}{Theorem}

\begingroup
\newtheorem{lemma}{Lemma}[section]
\newtheorem{proposition}[lemma]{Proposition}

\endgroup

\providecommand{\keywords}[1]{\textbf{Keywords: } #1}

\newtheorem{definition}{Definition}
{\theorembodyfont{\rmfamily}\newtheorem{remark}{Remark}[section]}

\newcommand{\ep}{\varepsilon}

\newcommand{\ds}{\displaystyle}
\newcommand{\beq}[1]{\begin{equation} \label{#1}\ds}
\newcommand{\eeq}{\end{equation}}
\newcommand{\bml}[1]{\beq{#1} \begin{array}{c}\ds}
\newcommand{\eml}{\end{array}\eeq}
\newcommand{\beqq}{\begin{equation*}\ds}
\newcommand{\eeqq}{\end{equation*}}
\newcommand{\bmll}{\beqq \begin{array}{c}\ds}
\newcommand{\emll}{\end{array}\eeqq}
\renewcommand{\div}{{\rm div}\,}

\newcommand{\abs}[1]{\ensuremath{\left| #1 \right|}}

\def \de{\partial}

\def \ep{\varepsilon}

\def \d{\mathrm{d}}

\newcommand{\R}{\mathbb{R}}

\newcommand{\id}{{\rm Id}}
\newcommand{\Sym}{\mathcal{S}}

\begin{document}

\author{Jan B\v rezina, Elisabetta Chiodaroli and Ond\v{r}ej Kreml\thanks{O.K. acknowledges the support of the GA\v CR (Czech Science Foundation) project GJ17-01694Y in the general framework of RVO: 67985840.}}
\title{On contact discontinuities in multi--dimensional isentropic Euler equations}
\date{}

\maketitle

\centerline{Tokyo Institute of Technology}

\centerline{2-12-1 Ookayama, Meguro-ku, Tokyo, 152-8550, Japan}

\bigskip

\centerline{EPFL Lausanne}

\centerline{Station 8, CH-1015 Lausanne, Switzerland}

\bigskip

\centerline{Institute of Mathematics, Czech Academy of Sciences}

\centerline{\v Zitn\'a 25, Prague 1, 115 67, Czech Republic}

 
\begin{abstract}
In this short note we partially extend the recent nonuniqueness results on admissible weak solutions to the Riemann problem for the 2D compressible isentropic Euler equations. We prove nonuniqueness of admissible weak solutions that start from the Riemann initial data allowing a contact discontinuity to emerge. 
\end{abstract}
\keywords{isentropic Euler equations; nonuniqueness; admissible weak solutions; Riemann problem; contact discontinuity}

\section{Introduction}\label{s:1}
Our concern in this paper is the isentropic compressible Euler system in two space dimensions
\begin{equation}\label{eq:Euler system}
\left\{\begin{array}{l}
\partial_t \rho + {\rm div}_x (\rho v) \;=\; 0\\
\partial_t (\rho v) + {\rm div}_x \left(\rho v\otimes v \right) + \nabla_x [ p(\rho)]\;=\; 0\\
\rho (\cdot,0)\;=\; \rho^0\\
v (\cdot, 0)\;=\; v^0. \, 
\end{array}\right.
\end{equation}
Here $(\rho,v)$ denote the unknown density and velocity of the fluid respectively. 
The pressure $p$ is a given function of $\rho$ and in order for  system \eqref{eq:Euler system} to be hyperbolic, it needs to satisfy $p'>0$.
Throughout this paper we assume that $p(\rho)=  \rho^\gamma$ with a constant $\gamma\geq 1$. The space variables are denoted as $x=(x_1, x_2)\in \R^2$ and similarly the components of the vectors are denoted as $v = (v_1,v_2) \in \R^2$.

The total energy of the fluid is given as the sum of the kinetic energy $\rho\frac{\abs{v}^2}{2}$ and the internal energy $\rho\ep(\rho)$ where the internal energy density $\ep(\rho)$ is related to the pressure through the relation $p(r)=r^2 \varepsilon'(r)$. The total energy plays the role of the (only one) mathematical entropy in the terminology of hyperbolic conservations laws, therefore we also consider the \textit{entropy (energy) inequality}
\begin{equation} \label{eq:energy inequality}
\de_t \left(\rho \varepsilon(\rho)+\rho
\frac{\abs{v}^2}{2}\right)+\div_x
\left[\left(\rho\varepsilon(\rho)+\rho
\frac{\abs{v}^2}{2}+p(\rho)\right) v \right]
\;\leq\; 0.
\end{equation}

In this note we work with bounded weak solutions that satisfy \eqref{eq:Euler system} in the sense of distributions.
Moreover, we say that a weak solution to \eqref{eq:Euler system} is \textit{admissible}, when it satisfies \eqref{eq:energy inequality} 
in the sense of distributions, more precisely we require the following inequality to hold for every nonnegative 
test function $\varphi\in C_c^{\infty}(\R^2\times [0,\infty))$:
\begin{align*} 
 &\int_0^\infty\int_{\R^2} \left[\left(\rho\varepsilon(\rho)+\rho \frac{\abs{v}^2}{2}\right)\de_t \varphi+\left(\rho\varepsilon(\rho)+\rho
\frac{\abs{v}^2}{2}+p(\rho)\right) v \cdot \nabla_x \varphi \right] \d x \d t\notag \\
&+\int_{\R^2} \left(\rho^0 (x) \varepsilon(\rho^0 (x))+\rho^0 (x)\frac{\abs{v^0 (x)}^2}{2}\right) 
\varphi(x,0) \d x \geq 0 .
\end{align*}  

Camillo De Lellis and L\'aszlo Sz\'ekelyhidi proved in \cite{dls2} the existence of initial data $(\rho^0, v^0)$ for which there exists infinitely many admissible weak solutions to \eqref{eq:Euler system} by a suitable application of their theory for the incompressible Euler equations based on convex integration or Baire cathegory method. Later in \cite{ch} and \cite{ChDLKr} the regularity of such initial data was improved. The proof in \cite{ChDLKr} uses as a core idea the analysis of the Riemann problem for compressible Euler equations in 2D. The Riemann problem is a problem with a specific choice of initial data in the following form
\begin{equation}\label{eq:R_data}
(\rho^0 (x), v^0 (x)) := \left\{
\begin{array}{ll}
(\rho_-, v_-) \quad & \mbox{if $x_2<0$}\\ \\
(\rho_+, v_+) & \mbox{if $x_2>0$,} 
\end{array}\right. 
\end{equation}
where $\rho_\pm, v_\pm$ are constants. The same problem was further studied in \cite{ChKr1} and \cite{ChKr2} and also by Klingenberg and Markfelder \cite{KlMa}. All these results show that the entropy inequality itself is not enough to single out a unique physical solution for certain ranges of the initial data $\rho_\pm, v_\pm$.

The Riemann problem is a classical building block of the one-dimensional theory for hyperbolic conservation laws. It is well known that it allows for existence of $BV$ self-similar solutions consisting of constant states joined by rarefaction waves, admissible shocks and contact discontinuities, see for example \cite{da}. Since the initial data \eqref{eq:R_data} are one-dimensional, it is easy to observe, that the 1D self-similar solutions prolonged as constant to  the next dimension are indeed solutions to the 2D problem as well. Such solutions are unique in the class of  admissible weak solutions if we require them to be self-similar and to have locally bounded variation. However dropping the requirements of self-similarity and $BV_{loc}$ yields nonuniqueness as was illustrated in \cite{ChDLKr}, \cite{ChKr1}, \cite{ChKr2}.

In the case of 2D isentropic Euler system, a contact discontinuity appears in the self-similar solution if and only if the first components of the velocities $v_-$ and $v_+$ are not equal. If $v_{-1} = v_{+1}$, then the self-similar solution can consist only of admissible shocks and rarefaction waves, see \cite[Section 2]{ChKr1} for detailed analysis.

If the self-similar solution consists only of rarefaction waves, it is in fact unique in the class of all bounded admissible weak solutions, as was first proved in \cite{chen} (see also \cite{fekr} and \cite{serre} for related results). If on the other hand the self-similar solution consists of two admissible shocks, then there exists infinitely many admissible weak solutions with the same initial data, see \cite{ChKr1}. The same nonuniqueness result holds also in some cases where the self-similar solution consists of one shock and one rarefaction wave, see \cite{ChKr2}.

The case of Riemann initial data including $v_{-1} \neq v_{+1}$ has not been studied in this context yet, even though there is an interesting result by Sz\'ekelyhidi \cite{sz} concerning  incompressible Euler system. He proved that vortex sheet initial data (i.e. $v_- = (-1,0)$, $v_+ = (1,0)$) allow for the existence of infinitely many weak solutions (to incompressible Euler system) satisfying either strict energy inequality or energy equality. However this result does not seem to transfer directly to the compressible case mainly because the role of the pressure is different in both systems of equations.

The question, whether a self-similar solution consisting only of a contact discontinuity (or more generally rarefaction waves and a contact discontinuity) is unique in the class of bounded admissible weak solutions or not, is to our knowledge still open. On the other hand it is natural to expect that the nonuniqueness results of \cite{ChKr1} and \cite{ChKr2} also extend to the case when the self-similar solution contains a contact discontinuity. In this note we give a confirmation of this conjecture and show how to obtain nonuniqueness of admissible weak solutions for such initial data.

Our main results are as follows.

\begin{theorem}\label{t:main}
Let $p(\rho) = \rho^\gamma$, $\gamma \geq 1$. Let $\rho_+,\rho_- > 0$, $v_+,v_- \in \mathbb{R}^2$ and let $v_{-2} - v_{+2} > \sqrt{\frac{(\rho_+ - \rho_-)(p(\rho_+) - p(\rho_-))}{\rho_+\rho_-}}$. Then there exists infinitely many admissible weak solutions to the Riemann problem for the Euler system \eqref{eq:Euler system} and \eqref{eq:R_data}.
\end{theorem}

\begin{remark}\label{r:1}
The condition $v_{-2} - v_{+2} > \sqrt{\frac{(\rho_+ - \rho_-)(p(\rho_+) - p(\rho_-))}{\rho_+\rho_-}}$ means that 
the self-similar solution consists of two shocks and a contact discontinuity (apart from the case $v_{-1} = v_{+1}$
 when the contact discontinuity does not appear). For a schematic view of the self-similar solution in this case see Figure \ref{fig:waves0}. 
Theorem \ref{t:main} extends the result of \cite{ChKr1} to the case $v_{-1} \neq v_{+1}$.
\begin{figure}[h] 
\begin{center}
\scalebox{0.5}{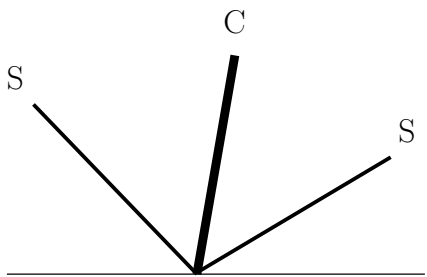}
\caption{Wave fan consisting of a shock, a contact discontinuity and another shock.}
\label{fig:waves0}
\end{center}
\end{figure}
\end{remark}

\begin{theorem}\label{t:main1}
Let $p(\rho) = \rho^\gamma$, $\gamma \geq 1$. Let $\rho_+,\rho_- > 0$, $\rho_+\neq\rho_-$, $v_+,v_- \in \mathbb{R}^2$ and     
  let $v_{+2} \in \mathbb{R}$. There exists $\overline{V} := \overline{V}(\rho_-,\rho_+,v_{+2},\gamma) < \sqrt{\frac{(\rho_+ - \rho_-)(p(\rho_+) - p(\rho_-))}{\rho_+\rho_-}}$ such that if 
\begin{equation*}
\overline{V} < v_{-2} - v_{+2} < \sqrt{\frac{(\rho_+ - \rho_-)(p(\rho_+) - p(\rho_-))}{\rho_+\rho_-}} 
\end{equation*}
then there exists infinitely many admissible weak solutions to the Riemann problem for the Euler system \eqref{eq:Euler system} and \eqref{eq:R_data}.
\end{theorem}

\begin{remark}\label{r:2}
The condition $\overline{V} < v_{-2} - v_{+2} < \sqrt{\frac{(\rho_+ - \rho_-)(p(\rho_+) - p(\rho_-))}{\rho_+\rho_-}}$ means that the self-similar solution consists of one shock, 
one rarefaction wave and a contact discontinuity (apart from the case $v_{-1} = v_{+1}$ when the contact discontinuity does not appear). See Figure \ref{fig:waves} for a 
sketch of the two possible structures of the self-similar solution in the case $v_{-1} \neq v_{+1}$. Theorem \ref{t:main1} extends the result 
of \cite{ChKr2} to the case $v_{-1} \neq v_{+1}$.
 \begin{figure}[hp]
 \centering
 \subfigure[]
  {\scalebox{0.5}{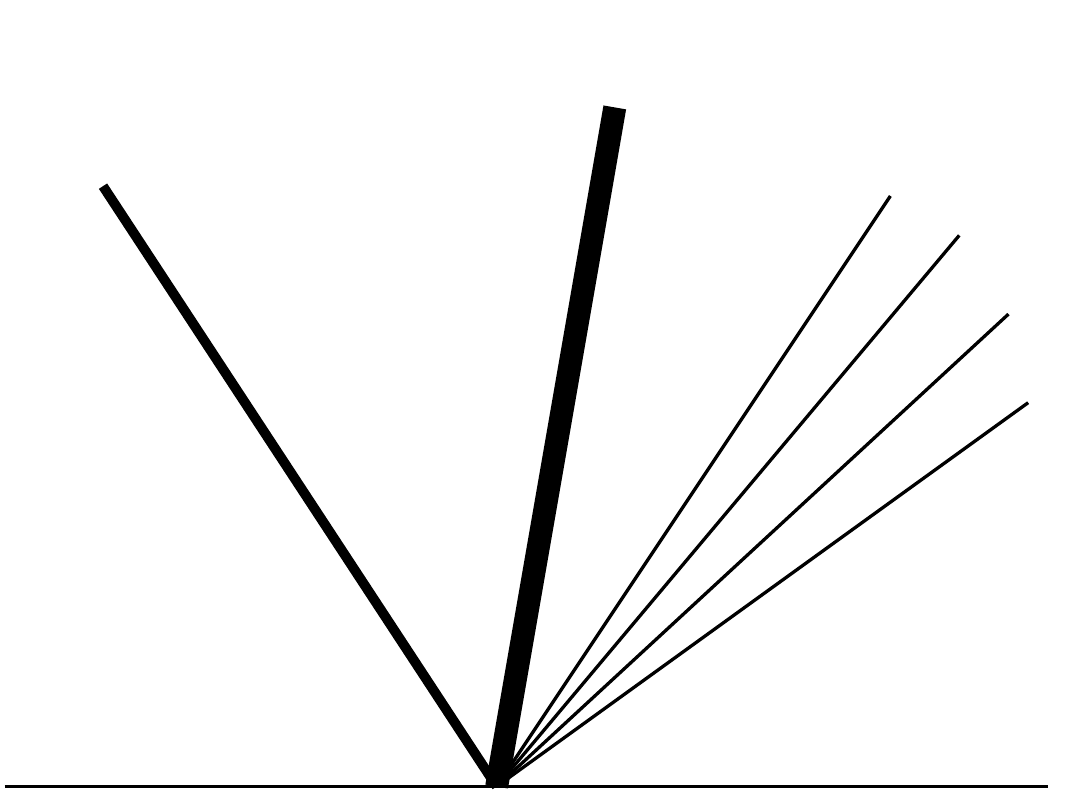}}
 \hspace{8mm}
 \subfigure[]
 {\scalebox{0.5}{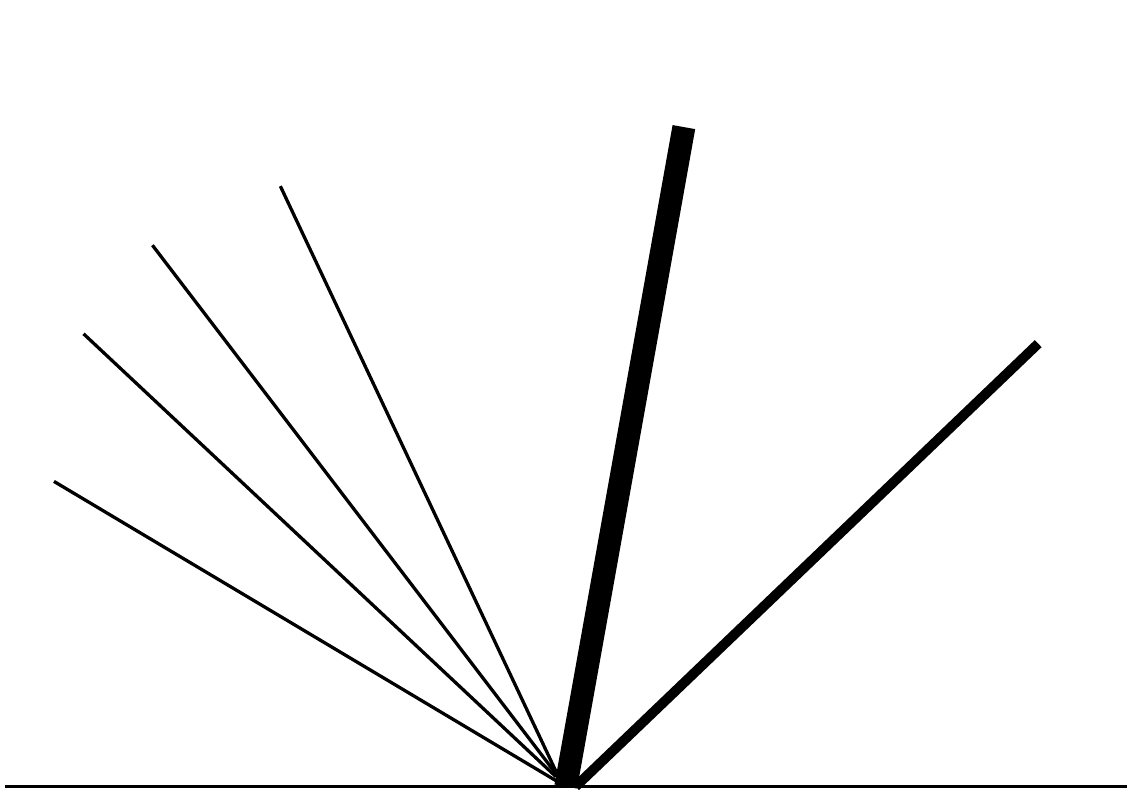}}
 \caption{Wave fan consisting of (a) a shock, a contact discontinuity and a rarefaction, or (b) a rarefaction, a contact discontinuity and a shock.}
\label{fig:waves}
 \end{figure}
\end{remark}

In the rest of this note we prove Theorem \ref{t:main} and Theorem \ref{t:main1}.

\section{Preliminaries}\label{s:2}

Here we state three important definitions from \cite{ChDLKr} in the form we need in this paper.

\begin{definition}[Fan partition]\label{d:fan}
A {\em fan partition} of $\R^2\times (0, \infty)$ consists of four open sets $P_-, P_1, P_2, P_+$
of the following form 
\begin{align*}
 P_- &= \{(x,t): t>0 \quad \mbox{and} \quad x_2 < \nu_- t\}\\
 P_1 &= \{(x,t): t>0 \quad \mbox{and} \quad \nu_- t < x_2 < \nu_1 t\}\\
 P_2&= \{(x,t): t>0 \quad \mbox{and} \quad \nu_1 t < x_2 < \nu_+ t\}\\
 P_+ &= \{(x,t): t>0 \quad \mbox{and} \quad x_2 > \nu_+ t\},
\end{align*}
where $\nu_- < \nu_1 < \nu_+$ is an arbitrary triple of real numbers.
\end{definition}

\begin{definition}[Fan subsolution] \label{d:subs}
A {\em fan subsolution} to the compressible Euler equations \eqref{eq:Euler system} with the 
initial data \eqref{eq:R_data} is a triple 
$(\overline{\rho}, \overline{v}, \overline{u}): \R^2\times 
(0,\infty) \rightarrow (\R^+, \R^2, \Sym_0^{2\times2})$ of piecewise constant functions satisfying
the following requirements.
\begin{itemize}
\item[(i)] There is a fan partition $P_-, P_1, P_2, P_+$ of $\R^2\times (0, \infty)$ such that
\[
(\overline{\rho}, \overline{v}, \overline{u})=  
(\rho_-, v_-, u_-) \bm{1}_{P_-}
+ (\rho_1, v_1, u_1) \bm{1}_{P_1} + (\rho_2, v_2, u_2) \bm{1}_{P_2}
+ (\rho_+, v_+, u_+) \bm{1}_{P_+}
\]
where $\rho_i, v_i, u_i$ are constants with $\rho_i >0$ ($i = 1,2$) and $u_\pm =
v_\pm\otimes v_\pm - \textstyle{\frac{1}{2}} |v_\pm|^2 \id$;
\item[(ii)] There exist positive constants $C_1, C_2$ such that
\begin{equation*} \label{eq:subsolution 2}
v_i\otimes v_i - u_i < \frac{C_i}{2} \id\, 
\end{equation*}
for $i = 1,2$;
\item[(iii)] The triple $(\overline{\rho}, \overline{v}, \overline{u})$ solves the following system in the
sense of distributions:
\begin{align}
&\partial_t \overline{\rho} + {\rm div}_x (\overline{\rho} \, \overline{v}) \;=\; 0\label{eq:continuity}\\
&\partial_t (\overline{\rho} \, \overline{v})+{\rm div}_x \left(\overline{\rho} \, \overline{u} 
\right) + \nabla_x \left( p(\overline{\rho})+\frac{1}{2}\left( \overline{\rho} |\overline{v}|^2 \bm{1}_{P_+\cup P_-} + \sum_{i=1}^2 C_i \rho_i
\bm{1}_{P_i}\right)\right)= 0.\label{eq:momentum}
\end{align}
\end{itemize}

\end{definition}

\begin{definition}[Admissible fan subsolution]\label{d:admiss}
 A fan subsolution $(\overline{\rho}, \overline{v}, \overline{u})$ is said to be {\em admissible}
if it satisfies the following inequality in the sense of distributions
\begin{align} 
&\de_t \left(\overline{\rho} \varepsilon(\overline{\rho})\right)+\div_x
\left[\left(\overline{\rho}\varepsilon(\overline{\rho})+p(\overline{\rho})\right) \overline{v}\right]
 + \de_t \left( \overline{\rho} \frac{|\overline{v}|^2}{2} \bm{1}_{P_+\cup P_-} \right)
+ \div_x \left(\overline{\rho} \frac{|\overline{v}|^2}{2} \overline{v} \bm{1}_{P_+\cup P_-}\right)\nonumber\\
&\qquad\qquad+ \sum_{i = 1}^2\left[\de_t\left(\rho_i \, \frac{C_i}{2} \, \bm{1}_{P_i}\right) 
+ \div_x\left(\rho_i \, \overline{v} \, \frac{C_i}{2}  \, \bm{1}_{P_i}\right)\right]
\;\leq\; 0\, .\label{eq:admissible subsolution}
\end{align}
\end{definition}

A sufficient condition for the existence of infinitely many admissible weak solutions is the existence of a single admissible fan subsolution as is stated in the following proposition.

\begin{proposition}\label{p:subs}
Let $p$ be any $C^1$ function and $(\rho_\pm, v_\pm)$ be such that there exists at least one
admissible fan subsolution $(\overline{\rho}, \overline{v}, \overline{u})$ of \eqref{eq:Euler system}
with the initial data \eqref{eq:R_data}. Then there are infinitely 
many bounded admissible weak solutions $(\rho, v)$ to \eqref{eq:Euler system} and \eqref{eq:R_data} such that 
$\rho=\overline{\rho}$ and $\abs{v}^2\bm{1}_{P_i} = C_i$ ($i = 1,2$).
\end{proposition}

The core of the proof of Proposition \ref{p:subs} is the following fundamental lemma.

\begin{lemma}\label{l:ci}
Let $(\tilde{v}, \tilde{u})\in \R^2\times \Sym_0^{2\times 2}$ and $C_0>0$ be such that $\tilde{v}\otimes \tilde{v}
- \tilde{u} < \frac{C_0}{2} \id$. For any open set $\Omega\subset \R^2\times \R$ there are infinitely many maps
$(\underline{v}, \underline{u}) \in L^\infty (\R^2\times \R , \R^2\times \Sym_0^{2\times 2})$ with the following property
\begin{itemize}
\item[(i)] $\underline{v}$ and $\underline{u}$ vanish identically outside $\Omega$;
\item[(ii)] $\div_x \underline{v} = 0$ and $\partial_t \underline{v} + \div_x \underline{u} = 0$;
\item[(iii)] $ (\tilde{v} + \underline{v})\otimes (\tilde{v} + \underline{v}) - (\tilde{u} + \underline{u}) = \frac{C_0}{2} \id$
a.e. on $\Omega$.
\end{itemize}
\end{lemma}
The proof of  Lemma \ref{l:ci} can be found in \cite[Section 4]{ChDLKr} and it is essentially based on the theory of De Lellis and Sz\'ekelyhidi in \cite{dls2} 
for the incompressible Euler system. We will not present the proof here.

Proposition \ref{p:subs} is proved using Lemma \ref{l:ci} in the following way. In each of the regions $P_1$, $P_2$ we use Lemma \ref{l:ci} with $(\tilde{v}, \tilde{u}) = (v_i,u_i)$ and $C_0 = C_i$ to obtain $\underline{v}_i$. Then it is not difficult to check that each couple $(\overline{\rho}, \overline{v} + \sum_{i=1}^2\underline{v}_i)$ is indeed an admissible weak solution to \eqref{eq:Euler system}. For a complete proof of Proposition \ref{p:subs}, we refer the reader to \cite[Section 3.3]{ChDLKr}.

\section{Proofs}\label{s:3}

\subsection{Proof of Theorem \ref{t:main}}

In order to prove Theorem \ref{t:main} it is enough to find a single admissible fan subsolution due to Proposition \ref{p:subs}. Therefore we introduce the following notation
\begin{align*} 
v_i &= (\alpha_i, \beta_i),\\
v_- &= (v_{-1}, v_{-2})\\
v_+ &= (v_{+1}, v_{+2})\\
u_i &=\left( \begin{array}{cc}
    \gamma_i & \delta_i \\
    \delta_i & -\gamma_i\\
    \end{array} \right)\, 
\end{align*}
for $i=1,2$. Since the fan subsolution is by definition formed by piecewise constant functions, the partial differential equations \eqref{eq:continuity}--\eqref{eq:admissible subsolution} transfer to a set of Rankin-Hugoniot conditions on each of the three interfaces of the fan partition. We have 
\begin{itemize}
\item Rankine-Hugoniot conditions on the left interface:
\begin{align}
&\nu_- (\rho_- - \rho_1) \, =\,  \rho_- v_{-2} -\rho_1  \beta_1 \label{eq:cont_left}  \\
&\nu_- (\rho_- v_{-1}- \rho_1 \alpha_1) \, = \, \rho_- v_{-1} v_{-2}- \rho_1 \delta_1  \label{eq:mom_1_left}\\
&\nu_- (\rho_- v_{-2}- \rho_1 \beta_1) \, = \,  
\rho_- v_{-2}^2 + \rho_1 \gamma_1 +p (\rho_-)-p (\rho_1) - \rho_1 \frac{C_1}{2}\, ;\label{eq:mom_2_left}
\end{align}
\item Rankine-Hugoniot conditions on the middle interface:
\begin{align}
&\nu_1 (\rho_1 - \rho_2) \, =\,  \rho_1  \beta_1 - \rho_2  \beta_2 \label{eq:cont_middle}  \\
&\nu_1 (\rho_1 \alpha_1 - \rho_2 \alpha_2) \, = \,  \rho_1 \delta_1 -  \rho_2 \delta_2  \label{eq:mom_1_middle}\\
&\nu_1 (\rho_1 \beta_1 - \rho_2\beta_2) \, = \,  
- \rho_1 \gamma_1 + \rho_2\gamma_2 + p (\rho_1)-p (\rho_2) + \rho_1 \frac{C_1}{2} - \rho_2 \frac{C_2}{2}\, ;\label{eq:mom_2_middle}
\end{align}
\item Rankine-Hugoniot conditions on the right interface:
\begin{align}
&\nu_+ (\rho_2-\rho_+ ) \, =\,  \rho_2  \beta_2 - \rho_+ v_{+2} \label{eq:cont_right}\\
&\nu_+ (\rho_2 \alpha_2 - \rho_+ v_{+1}) \, = \, \rho_2 \delta_2 - \rho_+ v_{+1} v_{+2} \label{eq:mom_1_right}\\
&\nu_+ (\rho_2 \beta_2 - \rho_+ v_{+2}) \, = \, - \rho_2 \gamma_2 - \rho_+ v_{+2}^2 +p (\rho_2) -p (\rho_+) 
+ \rho_2 \frac{C_2}{2}\, ;\label{eq:mom_2_right}
\end{align}
\item Subsolution conditions:
\begin{align}
&\alpha_1^2 +\beta_1^2 < C_1 \label{eq:sub_trace_1}\\
&\alpha_2^2 +\beta_2^2 < C_2 \label{eq:sub_trace_2}\\
& \left( \frac{C_1}{2} -{\alpha_1}^2 +\gamma_1 \right) \left( \frac{C_1}{2} -{\beta_1}^2 -\gamma_1 \right) - 
\left( \delta_1 - \alpha_1 \beta_1 \right)^2 >0\, \label{eq:sub_det_1}\\
& \left( \frac{C_2}{2} -{\alpha_2}^2 +\gamma_2 \right) \left( \frac{C_2}{2} -{\beta_2}^2 -\gamma_2 \right) - 
\left( \delta_2 - \alpha_2 \beta_2 \right)^2 >0\, ;\label{eq:sub_det_2}
\end{align}
\item Admissibility condition on the left interface:
\begin{align}\label{eq:E_left}
& \nu_-(\rho_- \varepsilon(\rho_-)- \rho_1 \varepsilon( \rho_1))+\nu_- 
\left(\rho_- \frac{\abs{v_-}^2}{2}- \rho_1 \frac{C_1}{2}\right)\\
\leq & \left[(\rho_- \varepsilon(\rho_-)+ p(\rho_-)) v_{-2}- 
( \rho_1 \varepsilon( \rho_1)+ p(\rho_1)) \beta_1 \right] 
+ \left( \rho_- v_{-2} \frac{\abs{v_-}^2}{2}- \rho_1 \beta_1 \frac{C_1}{2}\right)\, ;\nonumber
\end{align}
\item Admissibility condition on the middle interface:
\begin{align}
& \nu_1(\rho_1 \varepsilon(\rho_1)- \rho_2 \varepsilon( \rho_2))+\nu_1 
\left(\rho_1 \frac{C_1}{2}- \rho_2 \frac{C_2}{2}\right)\nonumber\\
\leq & \left[(\rho_1 \varepsilon(\rho_1)+ p(\rho_1)) \beta_1- 
( \rho_2 \varepsilon( \rho_2)+ p(\rho_2)) \beta_2 \right] 
+ \left( \rho_1 \beta_1 \frac{C_1}{2}- \rho_2 \beta_2 \frac{C_2}{2}\right)\, ;\label{eq:E_middle}
\end{align}
\item Admissibility condition on the right interface:
\begin{align}\label{eq:E_right}
&\nu_+(\rho_2 \varepsilon( \rho_2)- \rho_+ \varepsilon(\rho_+))+\nu_+ 
\left( \rho_2 \frac{C_2}{2}- \rho_+ \frac{\abs{v_+}^2}{2}\right)\\
\leq &\left[ ( \rho_2 \varepsilon( \rho_2)+ p(\rho_2)) \beta_2- (\rho_+ \varepsilon(\rho_+)+ p(\rho_+)) v_{+2}\right] 
+ \left( \rho_2 \beta_2 \frac{C_2}{2}- \rho_+ v_{+2} \frac{\abs{v_+}^2}{2}\right)\, .\nonumber
\end{align}
\end{itemize}

Motivated both by the structure of the self-similar solution as well as the structure of the fan subsolution from \cite{ChKr1} in the case of no contact discontinuity we make the following ansatz. We set
\begin{align}
&\alpha_1 = v_{-1}\label{eq:ansatz1} \\
&\alpha_2 = v_{+1} \\
&\rho_1 = \rho_2 \\
&\beta_1 = \beta_2 =: \beta.\label{eq:ansatz5}
\end{align}
Such ansatz yields the following simplification of the above set of identities and inequalities. The equation \eqref{eq:cont_middle} is satisfied trivially and the equation \eqref{eq:mom_2_middle} simplifies to 
\begin{equation} \label{eq:gamma12C1C2}
 \gamma_1 - \frac{C_1}2 = \gamma_2 - \frac{C_2}2.
\end{equation}
 Moreover combining \eqref{eq:cont_left} and \eqref{eq:mom_1_left} yields $\delta_1 = \alpha_1\beta$ and similarly we get from the right interface that $\delta_2 = \alpha_2\beta$. Plugging this in \eqref{eq:mom_1_middle} leads to $\nu_1 = \beta$. Finally the admissibility condition on the middle interface \eqref{eq:E_middle} is trivially satisfied. Thus, after applying \eqref{eq:gamma12C1C2} what remains is the following set of relations.
\begin{itemize}
\item Rankine-Hugoniot conditions on the left interface:
\begin{align}
&\nu_- (\rho_- - \rho_1) \, =\,  \rho_- v_{-2} -\rho_1  \beta \label{eq:2cont_left}  \\
&\nu_- (\rho_- v_{-2}- \rho_1 \beta) \, = \,  
\rho_- v_{-2}^2 + \rho_1 \gamma_1 +p (\rho_-)-p (\rho_1) - \rho_1 \frac{C_1}{2}\, ;\label{eq:2mom_2_left}
\end{align}
\item Rankine-Hugoniot conditions on the right interface:
\begin{align}
&\nu_+ (\rho_1-\rho_+ ) \, =\,  \rho_1  \beta - \rho_+ v_{+2} \label{eq:2cont_right}\\
&\nu_+ (\rho_1 \beta - \rho_+ v_{+2}) \, = \, - \rho_1 \gamma_1 - \rho_+ v_{+2}^2 +p (\rho_1) -p (\rho_+) 
+ \rho_1 \frac{C_1}{2}\, ;\label{eq:2mom_2_right}
\end{align}
\item Subsolution conditions:
\begin{align}
&v_{-1}^2 +\beta^2 < C_1 \label{eq:2sub_trace_1}\\
&v_{+1}^2 +\beta^2 < C_2 \label{eq:2sub_trace_2}\\
& \left( \frac{C_1}{2} -v_{-1}^2 +\gamma_1 \right) \left( \frac{C_1}{2} -{\beta}^2 -\gamma_1 \right) >0\, \label{eq:2sub_det_1}\\
& \left( \frac{C_2}{2} -v_{+1}^2 +\gamma_2 \right) \left( \frac{C_1}{2} -{\beta}^2 -\gamma_1 \right) >0\, ;\label{eq:2sub_det_2}
\end{align}
\item Admissibility condition on the left interface:
\begin{align}\label{eq:2E_left}
& \nu_-(\rho_- \varepsilon(\rho_-)- \rho_1 \varepsilon( \rho_1))+\nu_- 
\left(\rho_- \frac{v_{-1}^2+v_{-2}^2}{2}- \rho_1 \frac{C_1}{2}\right)\\
\leq & \left[(\rho_- \varepsilon(\rho_-)+ p(\rho_-)) v_{-2}- 
( \rho_1 \varepsilon( \rho_1)+ p(\rho_1)) \beta \right] 
+ \left( \rho_- v_{-2} \frac{v_{-1}^2+v_{-2}^2}{2}- \rho_1 \beta \frac{C_1}{2}\right)\, ;\nonumber
\end{align}
\item Admissibility condition on the right interface:
\begin{align}\label{eq:2E_right}
&\nu_+(\rho_1 \varepsilon( \rho_1)- \rho_+ \varepsilon(\rho_+))+\nu_+ 
\left( \rho_1 \frac{C_2}{2}- \rho_+ \frac{v_{+1}^2+v_{+2}^2}{2}\right)\\
\leq &\left[ ( \rho_1 \varepsilon( \rho_1)+ p(\rho_1)) \beta- (\rho_+ \varepsilon(\rho_+)+ p(\rho_+)) v_{+2}\right] 
+ \left( \rho_1 \beta \frac{C_2}{2}- \rho_+ v_{+2} \frac{v_{+1}^2+v_{+2}^2}{2}\right)\, .\nonumber
\end{align}
\end{itemize}

As argued in \cite[Lemma 4.3]{ChKr1} the inequalities \eqref{eq:2sub_trace_1}--\eqref{eq:2sub_det_2} are satisfied only if $$\frac{C_1}2 - \gamma_1 >  \beta^2.$$ Hence using the notation 
\begin{align*}
&\varepsilon_1 = \frac{C_1}2 - \gamma_1 - \beta^2\\
&\varepsilon_2 = \frac{C_1}2 - v_{-1}^2 +\gamma_1  = C_1  - v_{-1}^2 - \beta^2 - \varepsilon_1\\
& \varepsilon_2' = \frac{C_2}2 - v_{+1}^2 +\gamma_2  = C_2  - v_{+1}^2 - \beta^2 - \varepsilon_1
\end{align*}
we see that \eqref{eq:2sub_trace_1}--\eqref{eq:2sub_det_2} are equivalent to
$\varepsilon_1 >0$, $\varepsilon_2 >0$ and $\varepsilon_2' >0$. Before we proceed any further let us set $\varepsilon_2 = \varepsilon_2'$, i.e. 
\begin{equation}\label{eq:C1v1C2v1} C_1 - v_{-1}^2 = C_2 - v_{+1}^2 .\end{equation}
Finally, following the proof of \cite[Lemma 4.4]{ChKr1} we rewrite \eqref{eq:2cont_left}--\eqref{eq:2E_right} in the new variables $\varepsilon_1$ and $\varepsilon_2$ as follows. 

\begin{itemize}
\item Rankine-Hugoniot conditions on the left interface:
\begin{align}
&\nu_- (\rho_- - \rho_1) \, =\,  \rho_- v_{-2} -\rho_1  \beta \label{eq:3cont_left}  \\
&\nu_- (\rho_- v_{-2}- \rho_1 \beta) \, = \,  
\rho_- v_{-2}^2 - \rho_1(\beta^2 + \varepsilon_1)+p (\rho_-)-p (\rho_1) \, ;\label{eq:3mom_2_left}
\end{align}
\item Rankine-Hugoniot conditions on the right interface:
\begin{align}
&\nu_+ (\rho_1-\rho_+ ) \, =\,  \rho_1  \beta - \rho_+ v_{+2} \label{eq:3cont_right}\\
&\nu_+ (\rho_1 \beta - \rho_+ v_{+2}) \, = \,  \rho_1 (\beta^2 + \varepsilon_1) - \rho_+ v_{+2}^2 +p (\rho_1) -p (\rho_+) 
\, ;\label{eq:3mom_2_right}
\end{align}
\item Subsolution conditions:
\begin{align}
&\varepsilon_1 > 0 \label{eq:3sub_trace_1}\\
& \varepsilon_2 > 0\, ; \label{eq:3sub_trace_2}
\end{align}
\item Admissibility condition on the left interface:
\begin{align}\label{eq:3E_left}
& (\beta - v_{-2}) \left(p(\rho_-) + p (\rho_1) -2 \rho_- \rho_1 \frac{\varepsilon(\rho_-) - \varepsilon(\rho_1)}{\rho_- - \rho_1}\right) \\
\leq & \varepsilon_1 \rho_1(v_{-2} + \beta) - (\varepsilon_1 + \varepsilon_2) \frac{\rho_- \rho_1 (\beta - v_{-2})}{\rho_- - \rho_1}\, ;\nonumber
\end{align}
\item Admissibility condition on the right interface:
\begin{align}\label{eq:3E_right}
& (v_{+2} - \beta ) \left(p(\rho_1) + p (\rho_+) -2 \rho_1 \rho_+ \frac{\varepsilon(\rho_1) - \varepsilon(\rho_+)}{\rho_1 - \rho_+}\right) \\
\leq & \varepsilon_1 \rho_1(v_{+2} + \beta) - (\varepsilon_1 + \varepsilon_2) \frac{\rho_1 \rho_+ (v_{+2} - \beta )}{\rho_1 - \rho_+}\, .\nonumber
\end{align}

\end{itemize}

Now it is easy to  observe that the set of relations  \eqref{eq:3cont_left}--\eqref{eq:3E_right} is   exactly the same as the set of relations (4.26)--(4.33) in \cite{ChKr1}. The existence of a solution to this set of relations in the case $v_{-2} - v_{+2} > \sqrt{\frac{(\rho_+ - \rho_-)(p(\rho_+) - p(\rho_-))}{\rho_+\rho_-}}$ is proved in \cite[Section 4]{ChKr1}. To conclude that this solution together with the ansatz \eqref{eq:ansatz1}-\eqref{eq:ansatz5} and \eqref{eq:C1v1C2v1} defines in fact an admissible fan subsolution in the sense of Definition \ref{d:subs}, we only have to verify that $\nu_- < \nu_1 = \beta < \nu_+$. Indeed, from \eqref{eq:2cont_left} and \eqref{eq:2cont_right} we deduce that
\begin{align*}
 \beta - \nu_- &= \frac{\rho_-}{\rho_1}\left(v_{-2}-\nu_-\right) \\
 \nu_+ - \beta &= \frac{\rho_+}{\rho_1}\left(\nu_+ - v_{+2}\right) 
\end{align*} 
and the proof is finished using \cite[Lemma 4.6]{ChKr1} which states that 
\begin{align*}
  v_{-2}-\nu_- &> 0 \\
  \nu_+ - v_{+2} &>0. 
 \end{align*}
This concludes the proof of Theorem \ref{t:main}.

\begin{remark}
It is not difficult to observe that \cite[Theorem 2]{ChKr1} transfers to our case as well and we obtain in particular that there exists a Riemann initial data \eqref{eq:R_data} with $v_{-1} \neq v_{+1}$ such that the self-similar solution to the Euler system \eqref{eq:Euler system}, \eqref{eq:energy inequality} is not entropy rate admissible. For the definition of entropy rate admissibility see \cite[Definition 1]{ChKr1}.
\end{remark}

\subsection{Proof of Theorem \ref{t:main1}}

Let us first recall  \cite[Theorem 1]{ChKr2} here.
\begin{theorem}\label{t:prev}
Let $p(\rho) = \rho^\gamma$, $\gamma > 1$. Let $\rho_- \neq \rho_+$, $\rho_{\pm} > 0$ and $v_{+2} \in \mathbb{R}$ be given and let $v_{-1} = v_{+1}$. 
There exists $V = V(\rho_-,\rho_+,v_{+2},\gamma) < \sqrt{\frac{(\rho_+ - \rho_-)(p(\rho_+) - p(\rho_-))}{\rho_+\rho_-}}$ 
such that for all $v_{-2}$ satisfying $V < v_{-2} - v_{+2} < \sqrt{\frac{(\rho_+ - \rho_-)(p(\rho_+) - p(\rho_-))}{\rho_+\rho_-}}$ there exists infinitely many bounded admissible weak solutions to the Euler equations \eqref{eq:Euler system} with Riemann initial data \eqref{eq:R_data}.
\end{theorem}
The proof in \cite{ChKr2} is based on the analysis of the set of identities and inequalities \eqref{eq:2cont_left}-\eqref{eq:2E_right} with the specific choice $v_{-1} = v_{+1}$. A solution is proved to exist under the condition in Theorem \ref{t:prev}. In order to prove Theorem \ref{t:main1} we again search for a single admissible fan subsolution and use the same ansatz \eqref{eq:ansatz1}-\eqref{eq:ansatz5} and \eqref{eq:C1v1C2v1} as in the proof of Theorem \ref{t:main}. We argue the same way as before and the only condition we have to ensure is that $\nu_- < \beta < \nu_+$. As it is described in \cite[Section 3]{ChKr2}, this is indeed the case at least on a small neighborhood of $\sqrt{\frac{(\rho_+ - \rho_-)(p(\rho_+) - p(\rho_-))}{\rho_+\rho_-}}$. However, as it is shown in the examples in \cite[Section 4]{ChKr2}, there are subsolutions violating the condition $\beta < \nu_+$, so requiring this to hold yields a more restrictive lower bound on $v_{-2} - v_{+2}$, i.e. in general $V(\rho_-,\rho_+,v_{+2},\gamma) \leq \overline{V}(\rho_-,\rho_+,v_{+2},\gamma)$. Finally, it is not difficult to observe by going through the proof in \cite{ChKr2} that Theorem \ref{t:prev} holds also for $\gamma = 1$.

Theorem \ref{t:main1} is proved.


\begin{thebibliography}{plain}


 


 
 
 

\bibitem{chen} CHEN, G.-Q., CHEN, J.: Stability of rarefaction waves and vacuum states for the multidimensional Euler equations. 
\textit{J. Hyperbolic Differ. Equ.} \textbf{4}, 105--122, (2007).

%

\bibitem{ch} CHIODAROLI, E.: A counterexample to well-posedeness of entropy solutions to the compressible Euler system. \textit{J. Hyperbolic Differ. Equ.} \textbf{11}, 493--519, (2014).

\bibitem{ChDLKr} CHIODAROLI, E., DE LELLIS, C., KREML, O.: Global ill-posedness of the isentropic system of gas dynamics. \textit{Comm. Pure Appl. Math.} \textbf{68}, 1157--1190, (2015).


\bibitem{ChKr1} CHIODAROLI, E., KREML, O.: On the energy dissipation rate of solutions to the compressible isentropic Euler system. \textit{Arch. Rational Mech. Anal.} \textbf{214}, 1019--1049, (2014).

\bibitem{ChKr2} CHIODAROLI, E., KREML, O.: Non-uniqueness of admissible weak solutions to the Riemann problem for the isentropic Euler equations. Preprint (2017), arXiv:1704.01747. [math.AP].


\bibitem{da} DAFERMOS, C.M.: \textit{Hyperbolic conservation laws in continuum physics}, vol. \textbf{325} 
of Grundleheren der Mathematischen Wissenschaften [Fundamental Principles of Mathematical Sciences]. Third edition.
Springer, Berlin, (2010).


%


\bibitem{dls2} DE LELLIS, C., SZ\'EKELYHIDI, L.J.: On admissibility criteria for weak solutions 
of the Euler equations. \textit{Arch. Ration. Mech. Anal.} \textbf{195}, 225--260, (2010).




\bibitem{fekr} FEIREISL, E.,  KREML, O.: Uniqueness of rarefaction waves in multidimensional compressible Euler system. \textit{J. Hyperbolic Differ. Equ.} \textbf{12}, 489--499, (2015).

\bibitem{KlMa} KLINGENBERG, C.,  MARKFELDER, S.: The Riemann problem for the multidimensional isentropic system of gas dynamics is ill-posed if it contains a shock. Preprint. (2017).


\bibitem{serre} SERRE, D.: Long-time stability in systems of conservation laws, using relative entropy/energy. \textit{Arch. Ration. Mech. Anal.} \textbf{219},  679--699, (2016). 


%
%
%

\bibitem{sz} SZ\'EKELYHIDI, L.J.: Weak solutions to the incompressible Euler equations with 
vortex sheet initial data. \textit{C. R. Acad. Sci. Paris}, Ser.I \textbf{349}, 1063--1066, (2011).
 
\end{thebibliography}
\end{document}